\newcommand{\R}{I\!\!R}

\newcommand{\D}{\mathcal D}

\newcommand{\dd}{\mathrm d}
\newcommand{\Ker}{\mathrm{Ker}}
\newcommand{\Img}{\mathrm{Im}}
\newcommand{\oo}{{\mathrm o}}

\makeatletter
\newcommand{\mysubsection}{\@startsection
{subsection}%
{2}%
{0mm}%
{-\baselineskip}%
{0.5\baselineskip}%
{\normalfont\normalsize\itshape}
}%
\makeatother

\documentclass[oneside,letterpaper,draft,10pt]{amsart}

\usepackage{amsmath}               
\usepackage{amsthm}                
\usepackage{times}
\usepackage{amscd}

\numberwithin{equation}{section}

\title[{\tiny Constrained variational problems}]%
{Lagrangian and Hamiltonian Formalism for Constrained Variational Problems}
\author[P.\ Piccione]{Paolo Piccione}
\address{Departamento de Matem\'atica,\hfill\break\indent  Universidade de S\~ao Paulo, Brazil}
\email{piccione@ime.usp.br}
\urladdr{http://www.ime.usp.br/\~{}piccione}
\author[D.\ Tausk]{Daniel V.\ Tausk}
\address{Departamento de Matem\'atica,\hfill\break\indent  Universidade de S\~ao Paulo, Brazil}
\email{tausk@ime.usp.br}
\urladdr{http://www.ime.usp.br/\~{}tausk}
\thanks{The authors are partially sponsored by CNPq, Brazil.}

\date{September 2001}

\begin{document}


\theoremstyle{plain}\newtheorem{teo}{Theorem}[section]
\theoremstyle{plain}\newtheorem{prop}[teo]{Proposition}
\theoremstyle{plain}\newtheorem{lem}[teo]{Lemma}
\theoremstyle{plain}\newtheorem{cor}[teo]{Corollary}
\theoremstyle{definition}\newtheorem{defin}[teo]{Definition}
\theoremstyle{remark}\newtheorem{rem}[teo]{Remark}
\theoremstyle{plain} \newtheorem{assum}[teo]{Assumption}
\theoremstyle{definition}\newtheorem{example}[teo]{Example}


\begin{abstract}
We consider solutions of Lagrangian variational problems with
linear constraints on the derivative. More precisely, given a
smooth distribution $\D\subset TM$ on $M$ and a time-dependent
Lagrangian $L$ defined on $\D$, we consider an action functional
$\mathcal L$ defined on the set $\Omega_{PQ}(M,\D)$ of horizontal
curves in $M$ connecting two fixed submanifolds $P,Q\subset M$.
Under suitable assumptions, the set $\Omega_{PQ}(M,\D)$ has the
structure of a smooth Banach manifold and we can thus study the
critical points of $\mathcal L$. If the Lagrangian $L$ satisfies
an appropriate {\em hyper-regularity\/} condition, we associate to
it a degenerate Hamiltonian $H$ on $TM^*$ using a general notion
of {\em Legendre transform\/} for maps on vector bundles. We prove
that the solutions of the Hamilton equations of $H$ are precisely
the critical points of $\mathcal L$. In the particular case where
$L$ is given by the quadratic form corresponding to a positive
definite metric on $\D$, we obtain the well-known characterization
of the {\em normal geodesics\/} in sub-Riemannian geometry (see
\cite{LS}); by adding a potential energy term to $L$, we reobtain
the equations of motion for the {\em Vakonomic mechanics\/} with
non holonomic constraints (see \cite{KO}).
\end{abstract}

\maketitle

\begin{section}{Introduction}
\label{sec:intro}

The aim of this paper is to generalize to the context of
constrained variational problems some classical results about the
correspondence between Lagrangian and Hamiltonian formalisms (see
for instance \cite{AM}). Particular cases of this theory are the
{\em sub-Riemannian geodesic problem} (see for instance \cite{LS,
Mo1, Mo3, PTJGP}), and the so called {\em Vakonomic\/} approach to
the non holonomic mechanics (see for instance \cite{AKN, CF, KO,
VG}).

The constrained variational problem studied is modelled by the
following setup: we consider an $n$-dimensional differentiable
manifold $M$ endowed with a smooth distribution $\D\subset TM$ of
rank $k$; moreover, we assume that it is given a (possibly
time-dependent) Lagrangian function $L$ on $\D$. In the non
holonomic mechanics, $M$ represents the configuration space, $\D$
the constraint, and $L$ is typically the difference between the
kinetic and a potential energy. In the sub-Riemannian geodesic
problem, $L$ is simply the quadratic form corresponding to a
positive definite metric on $\D$.

The solutions of the constrained variational problem are given by
curves $\gamma:[a,b]\to M$ that are critical points of the action
functional $\mathcal
L(\gamma)=\int_a^bL\big(t,\gamma(t),\dot\gamma(t)\big)\,\mathrm
dt$ defined on the space:
\[\Omega_{PQ}\big([a,b],M,\D\big)=\big\{\gamma:[a,b]\xrightarrow{\;\;C^1\;}M:\gamma(a)\in P,\
\gamma(b)\in Q,\ \gamma'(t)\in\D\ \text{for all}\ t\big\}\] of
horizontal curves of class $C^1$ in $M$ connecting two fixed
submanifolds $P,Q\subset M$. It is well-known that the set
$\Omega_{PQ}\big([a,b],M,\D\big)$ is in general not a submanifold
of the Banach manifold of $C^1$ curves $\gamma:[a,b]\to M$; when
$P$ and $Q$ are points, the singularities of $\Omega_{PQ}\big([a,b],M,\D\big)$ are known in the
context of sub-Riemannian geometry as {\em abnormals extremals\/}
(see \cite{BH, LS, Mo1, Mo2, Mo3, PTJGP}). Such singularities can
be nicely described using the canonical symplectic structure of the
cotangent bundle $TM^*$ (see Corollary~\ref{thm:outrosing}). In this paper we are interested in studying
the action functional $\mathcal L$ in the {\em regular part\/} of
$\Omega_{PQ}\big([a,b],M,\D\big)$. We remark that in several
important cases the set $\Omega_{PQ}\big([a,b],M,\D\big)$ contains
no singular curves (see, for instance, Corollary~\ref{thm:nosingular} and
Remark~\ref{thm:stronglybracket}).

Recall from \cite{AM} that when a Lagrangian function $L:TM\to\R$
is {\em hyper-regular\/} then the critical points of the
corresponding (unconstrained) variational problem are given by the
solutions of the Hamilton equations corresponding to a Hamiltonian
$H:TM^*\to\R$ which corresponds to $L$ by means of the {\em
Legendre transform}. The Legendre transform described in \cite{AM}
can be generalized in a straightforward way to general vector
bundles; namely, if $L:\xi\to\R$ is a smooth map on a vector
bundle $\xi$ which is hyper-regular (in a suitable sense) then one
can naturally associate to it a smooth map $H:\xi^*\to\R$ on the
dual bundle $\xi^*$. At such level of generality, the Legendre
transform does not seem to have a meaningful interpretation in the
context of calculus of variations, as it does in the case
$\xi=TM$. Our goal is to show that when $\xi=\D$ is a vector
subbundle of a tangent bundle $TM$ (i.e., a distribution on $M$)
then the Legendre transform for smooth maps on $\D$ has a nice
application to the study of constrained variational problems. The
key observation here is that, when passing to the dual bundles,
the {\em inclusion\/} arrow $\D\to TM$ reverses and gives rise to
a {\em projection\/} arrow $TM^*\to\D^*$; thus, while a {\em
constrained Lagrangian\/} $L:\D\to\R$ has no canonical extension
to a Lagrangian on $TM$, its Legendre transform $H_0:\D^*\to\R$
naturally induces a map $H:TM^*\to\R$ given by the composition of
$H_0$ and the projection $TM^*\to\D^*$. Our main result
(Theorem~\ref{thm:degenerate}) is that the critical points of the
constrained action functional $\mathcal L$ are the solutions of
the Hamilton equations of $H$ satisfying suitable boundary
conditions. Observe that, unless $\D=TM$, the Hamiltonian $H$ is
{\em always degenerate\/} and thus it cannot arise as the Legendre
transform of a hyper-regular Lagrangian on the whole tangent
bundle $TM$.

In the particular case where $P$ and $Q$ are single points of $M$,
$\D$ is endowed with a smoothly varying positive definite inner
product $g$ and $L$ is given by $L(t,q,\dot q)=\frac12 g(\dot
q,\dot q)$, then the solutions of the corresponding Hamiltonian
$H$ are known in the context of {\em sub-Riemannian geometry\/} as
the {\em normal extremals\/} of $(M,\D,g)$. The critical points of
the constraint defining $\Omega_{PQ}\big([a,b],M,\D\big)$ are
the abnormal extremals. In particular, we obtain a
variational proof of \cite[Theorem~1]{LS}. By adding a potential
energy term to $L$, the Hamilton equations of $H$ become the
equations of motion for the Vakonomic mechanics (see
\cite{KO}). Theorem~\ref{thm:degenerate} thus provides a unifying
approach for the study of Lagrangian variational problems with
linear constraints in the derivative; it also provides the
appropriate setting for the study of the {\em second variation\/}
of a constrained Lagrangian action functional and for the
development of an {\em index theory\/} for such functional using
the notion of {\em Maslov index\/} for a solution of a Hamiltonian
(see \cite{PTJMPA}).

The proof of Theorem~\ref{thm:degenerate} is based on the method
of Lagrangian multipliers, which is used to pass from a
constrained Lagrangian variational problem to a non constrained
one. The main technical difficulty is the proof of the regularity
of the Lagrangian multiplier (Lemma~\ref{thm:postponed}); such
proof is based on a suitable version of Schwartz's
generalized functions calculus which is developed in
Subsection~\ref{sub:calculus}.

We give a brief description of the material presented in each
section of the paper.

In Subsection~\ref{sub:legendre} we describe a general notion of
Legendre transform. In Subsection~\ref{sub:timedep} we
recall some standard results concerning the correspondence
between hyper-regular Lagrangians and Hamiltonians
and in Section~\ref{sec:horizontal} we present some
well-known facts about the manifold structure of the set of
horizontal curves connecting two fixed submanifolds of a
given manifold. 

In Section~\ref{sec:degenerate} we state the
main result of the paper (Theorem~\ref{thm:degenerate}),
that establishes the correspondence between the critical
points of the action functional of a hyper-regular
constrained Lagrangian  and the solutions of the
corresponding degenerate Hamiltonian. The proof of
Theorem~\ref{thm:degenerate} is given in
Subsection~\ref{sub:proof}. In Subsection~\ref{sub:calculus}
it is presented a suitable version of Schwartz's generalized
functions calculus, needed for technical reasons in the
proof of Theorem~\ref{thm:degenerate}.
\end{section}

\begin{section}{The Legendre Transform.\\ Lagrangians and Hamiltonians on Manifolds}
\label{sec:legendre}

In this section we recall some classical results from \cite{AM}
which are presented in a more general context  needed for
the statement and the proof of Theorem~\ref{thm:degenerate}. In
Subsection~\ref{sub:legendre} we present a general version of the
Legendre transform for vector spaces; we then apply it fiberwise
to obtain a notion of Legendre transform for fiber bundles. In
Subsection~\ref{sub:timedep} we present the classical Hamiltonian
formulation for the variational problem corresponding to a
hyper-regular (non constrained) Lagrangian. The standard results
from \cite{AM} are proven in a slightly more general
setup; namely, we consider curves with endpoints
varying in submanifolds, time-dependent Lagrangians and
rather weak regularity assumptions for the data.

\mysubsection{The Legendre transform}\label{sub:legendre} Let
$\xi_0$ be a real finite-dimensional  vector space, let $\xi_0^*$
denote its dual,  and let $Z:U\to\R$ be a function of class $C^2$
defined on an open subset $U\subset\xi_0$.
\begin{defin}\label{thm:defleg}
Assume that the differential $\mathrm dZ$ is a diffeomorphism onto
an open subset $V\subset\xi_0^*$. The {\em Legendre transform\/}
of $Z$ is the $C^1$ map $Z^*:V\to\R$ defined by:
\begin{equation}\label{eq:defZ*}
Z^*=E_Z\circ(\mathrm dZ)^{-1},
\end{equation}
where $E_Z:U\to\R $ is given by
\begin{equation}\label{eq:defEZ}
E_Z(v)=\mathrm dZ(v)\,v-Z(v),\quad v\in U.
\end{equation}
\end{defin}
\begin{lem}\label{thm:Z*C2}
Using the canonical identification of $\xi_0$ and its bi-dual $\xi_0^{**}$,
the map $\mathrm dZ^*$ is the inverse of $\mathrm dZ$. Therefore, $Z^*$ is a
map of class $C^2$.
\end{lem}
\begin{proof}
Differentiating the equality $Z^*\circ\mathrm dZ=E_Z$ and \eqref{eq:defEZ}, we obtain:
\[\mathrm dZ^*\big(\mathrm dZ(v)\big)\circ \mathrm d^2Z(v)=\mathrm dE_Z(v),\quad
\mathrm dE_Z(v)=\hat v\circ \mathrm d^2Z(v),\] where $\hat
v\in\xi_0^{**}$ denotes evaluation at $v$. Since $\mathrm
d^2Z(v):\xi_0\to\xi_0^*$ is an isomorphism, the conclusion
follows.
\end{proof}
\begin{cor}\label{thm:Z**Z}
$Z^{**}=Z$.
\end{cor}
\begin{proof}
By Lemma~\ref{thm:Z*C2}, we have:
\[Z^{**}=E_{Z^*}\circ(\mathrm dZ^*)^{-1}=E_{Z^*}\circ\mathrm dZ.\]
Hence, by definition of $E_{Z^*}$, we get:
\[\begin{split}
E_{Z^*}\big(\mathrm dZ(v)\big)\,&=\mathrm dZ^*\big(\mathrm
dZ(v)\big)\,\mathrm dZ(v)-Z^*\big(\mathrm dZ(v)\big)=\\&=\mathrm
dZ(v)\,v-E_Z(v)=Z(v).\qedhere
\end{split}\]
\end{proof}

Let now $M$ be a smooth manifold and $\pi:\xi\to M$ be a smooth
vector bundle over $M$; for $m\in M$, we denote by $\xi_m$ the
fiber $\pi^{-1}(m)$. The dual bundle of $\xi$ will be denoted by
$\xi^*$; the bi-dual $\xi^{**}$ is canonically identified with
$\xi$.

Let $Z:U\subset\xi\to \R$ be a map such that, for every $m\in M$,
$U\cap\xi_m$ is open in $\xi_m$ and the restriction of $Z$ to
$U\cap\xi_m$ is of class $C^2$.
\begin{defin}\label{thm:defFZ}
The {\em fiber derivative\/} $\mathbb FZ:U\to\xi^*$ is the map
defined by:
\begin{equation}\label{eq:defFZ}
\mathbb FZ(v)=\mathrm d(Z\vert_{U\cap\xi_m})(v),\quad v\in U,
\end{equation}
where $m=\pi(v)$. Let $V\subset\xi^*$ denote the image of $\mathbb
FZ$. We say that $Z$ is {\em regular\/} if for each $m\in M$, the
set $V\cap\xi_m$ is open in $\xi_m$ and the restriction of
$\mathbb FZ$ to $U\cap\xi_m$ is a local diffeomorphism; $Z$ is
said to be {\em hyper-regular\/} if for each $m$ such restriction
is a diffeomorphism onto $V\cap\xi^*_m$. If $Z$ is hyper-regular,
we define the {\em Legendre transform\/} of $Z$ as the map
$Z^*:V\to\R$ whose restriction to $V\cap\xi_m$ is the Legendre
transform of the restriction of $Z$ to $U\cap\xi_m$.
\end{defin}
In analogy with
\eqref{eq:defEZ} we also set:
\begin{equation}\label{eq:EZfibrado}
E_Z(v)=\mathbb FZ(v)\,v-Z(v),\quad v\in U;
\end{equation}
obviously $Z^*=E_Z\circ\mathbb FZ^{-1}$.

Applying Lemma~\ref{thm:Z*C2} and Corollary~\ref{thm:Z**Z}
fiberwise, we obtain immediately the following:
\begin{prop}\label{thm:bundles}
Assume that $Z:U\subset\xi\to\R$ is hyper-regular. Then, for each
$m\in M$, the restriction of $Z^*$ to $V\cap\xi^*_m$ is of class
$C^2$. Moreover, $\mathbb FZ$ and $\mathbb FZ^*$ are mutually
inverse bijections and $Z^{**}=Z$.\qed
\end{prop}

\mysubsection{Time dependent Lagrangians and Hamiltonians on
manifolds} \label{sub:timedep} Let $M$ be a smooth $n$-dimensional
manifold and let $TM$, $TM^*$ denote respectively the tangent and
the cotangent bundle of $M$; with a slight abuse of notation, we
will denote both the projections of $TM$ and of $TM^*$ by $\pi$.
Consider the following vector bundles:
\[
\xi=\R\times TM\xrightarrow{\;\;\mathrm{Id}\times\pi\;\;}\R\times
M,\quad\xi^*=\R\times
TM^*\xrightarrow{\;\;\mathrm{Id}\times\pi\;\;}\R\times M.\]
Observe that the fiber $\xi_{(t,m)}$ is $\{t\}\times T_mM$ and
that $\xi_{(t,m)}^*=\{t\}\times T_mM^*$.
\begin{defin}\label{thm:defLagrHam}
A {\em (time-dependent) Lagrangian on $M$\/} is a function
$L:U\subset\xi\to\R$ defined on an open set $U\subset\xi$ and
satisfying the following regularity conditions:
\begin{enumerate}
\item $L$ is continuous;
\item for each $t\in\R$, the map $L(t,\cdot)$ is of class $C^1$ on $U\cap\big(\{t\}\times TM\big)$
and its differential is continuous on $U$;
\item for each $t\in\R$, the map $\mathbb FL(t,\cdot):U\cap\big(\{t\}\times TM\big)
\to\{t\}\times TM^*$ is of class $C^1$.
\end{enumerate}
A {\em (time-dependent) Hamiltonian on $M$\/} is a function
$H:V\subset\xi^*\to\R$ defined on an open set $V\subset\xi^*$ and
satisfying the following regularity conditions:
\begin{enumerate}
\item for all $t\in\R$, the map $H(t,\cdot)$ is of class $C^1$ on $V\cap\big(\{t\}\times TM^*\big)$;
\item for each $(t,m)\in\R\times M$, the restriction of $H$ to $V\cap\xi^*_{(t,m)}$
is of class $C^2$.
\end{enumerate}
\end{defin}
We use the notions of regularity and hyper-regularity given in Definition~\ref{thm:defFZ}
for Lagrangians and Hamiltonians on manifolds.

Using the Legendre transform defined in
Subsection~\ref{sub:legendre} (Definition~\ref{thm:defFZ}), given
a hyper-regular Lagrangian $L$ on $M$, the map $H=L^*$ is a
hyper-regular Hamiltonian on $M$. Namely, the fact that
$H(t,\cdot)$ is of class $C^1$ follows by applying the Inverse
Function Theorem to the map $\mathbb FL(t,\cdot)$; moreover, the
fact that $V=\mathbb FL(U)$ is open in $\xi^*$ follows from the
Theorem of
Invariance of Domain (see \cite{Munkres}) by observing that $\mathbb FL$ is continuous and injective\footnote{%
As a matter of fact, this same argument shows that $\mathbb
FL:U\to V$ is a homeomorphism and therefore the Hamiltonian
$H=L^*$ is continuous.}.

If $H$ is the hyper-regular Hamiltonian obtained by Legendre
transform from the Lagrangian $L$, then by
Proposition~\ref{thm:bundles}, we have that $H^*=L$, and that
$\mathbb FH$ and $\mathbb FL$ are mutually inverse bijections. In
order to simplify the notation, in what follows we will write:
\[\mathbb FL(t,v)=\big(t,\mathbb FL^{(2)}(t,v)\big),\quad\mathbb
FH(t,p)=\big(t,\mathbb FH^{(2)}(t,p)\big),\] so that $\mathbb FL^{(2)}$
and $\mathbb FH^{(2)}$ are respectively a $TM^*$-valued and a
$TM$-valued map.

Let $L:U\subset\R\times TM\to\R$ be a Lagrangian on $M$ and
$\gamma:[a,b]\to M$ be a curve of class $C^1$, with
$\big(t,\dot\gamma(t)\big)\in U$ for all $t$. The {\em action\/}
$\mathcal L(\gamma)$ of $L$ on the curve $\gamma$ is given by the
integral:
\begin{equation}\label{eq:action}
\mathcal L(\gamma)=\int_a^bL\big(t,\dot\gamma(t)\big)\,\mathrm dt.
\end{equation}
$\mathcal L$ defines a functional on the set:
\begin{multline}\label{eq:set}
\Omega_{PQ}\big([a,b],M;U\big)\\
=\Big\{\gamma:[a,b]\xrightarrow{\;\;C^1\;}M:\gamma(a)\in P,\
\gamma(b)\in Q,\ \big(t,\dot\gamma(t)\big)\in U,\
\forall\,t\in[a,b]\Big\},
\end{multline}
where $P$ and $Q$ are two smooth embedded submanifolds of $M$. It
is well known that $\Omega_{PQ}\big([a,b],M;U\big)$ has the
structure of an infinite dimensional smooth Banach manifold (see
for instance \cite{Palais, PT1}), and $\mathcal L$ is a functional
of class $C^1$ on $\Omega_{PQ}\big([a,b],M;U\big)$. We will call
$\mathcal L$ the {\em action functional\/} associated to the
Lagrangian $L$.

We have the following characterization of the critical points of $\mathcal L$:
\begin{prop}\label{thm:EL}
A curve $\gamma\in\Omega_{PQ}\big([a,b],M;U\big)$ is a critical
point of $\mathcal L$ if and only if the following three
conditions are satisfied:
\begin{enumerate}
\item\label{itm:1} $\mathbb FL^{(2)}\big(a,\dot\gamma(a)\big)\vert_{T_{\gamma(a)P}}=0$ and $\mathbb
FL^{(2)}\big(b,\dot\gamma(b)\big)\vert_{T_{\gamma(b)Q}}=0$;
\item\label{itm:2} $t\mapsto \mathbb FL\big(t,\dot\gamma(t)\big)$ is of class $C^1$;
\item\label{itm:3} for all $[t_0,t_1]\subset[a,b]$ and for any chart
$q=(q_1,\ldots,q_n)$ on $M$ whose domain contains
$\gamma\big([t_0,t_1]\big)$, the {\em Euler--Lagrange equation\/}
is satisfied in $[t_0,t_1]$:
\begin{equation}\label{eq:EL}
\frac{\mathrm d}{\mathrm dt}\,\frac{\partial L}{\partial\dot
q}\big(t,q(t),\dot q(t)\big)= \frac{\partial L}{\partial
q}\big(t,q(t),\dot q(t)\big),
\end{equation}
where $L(t,q,\dot q)$ denotes the coordinate representation of
$L$.
\end{enumerate}
\end{prop}
\begin{proof}
Let $\gamma\in\Omega_{PQ}\big([a,b],M;U\big)$ be a critical point
of $\mathcal L$. Let $[t_0,t_1]\subset[a,b]$ be an interval and
consider a chart $q=(q_1,\ldots,q_n)$ in $M$ whose domain contains
$\gamma\big([t_0,t_1]\big)$. Choose an arbitrary $v\in
T_\gamma\Omega_{PQ}\big([a,b],M;U\big)$ with support contained in
$\left]t_0,t_1\right[$; by standard computations it follows that:
\begin{equation}\label{eq:standardcomp}
\int_{t_0}^{t_1}\frac{\partial L}{\partial q}\big(t,q(t),\dot
q(t)\big)\,v(t)+ \frac{\partial L}{\partial \dot
q}\big(t,q(t),\dot q(t)\big)\,\dot v(t)\;\mathrm dt=0.
\end{equation}
The fact that the equality above holds for every smooth $v$ with
support contained in $\left]t_0,t_1\right[$ implies that the term
$\frac{\partial L}{\partial \dot q}\big(t,q(t),\dot q(t)\big)$ is
of class $C^1$; this will follow\footnote{%
Alternatively, one could use integration by parts and the fact
that $\int_{t_0}^{t_1}\phi\dot v=0$ for all smooth $v$ with
support in $\left]t_0,t_1\right[$ implies
$\phi\equiv\text{constant}$.} from the generalized functions
calculus developed in Subsection~\ref{sub:calculus} (see
Corollary~\ref{thm:bootstrap}). Integration by parts in
\eqref{eq:standardcomp} and the Fundamental Lemma of Calculus of
Variations imply then that equation \eqref{eq:EL} is satisfied.
Observe also that the coordinate representation of the map
$t\mapsto\mathbb FL^{(2)}\big(t,\dot\gamma(t)\big)$ is given by
$t\mapsto\frac{\partial L}{\partial \dot q}\big(t,q(t),\dot
q(t)\big)$, so that condition~\eqref{itm:2} is satisfied.
Condition~\eqref{itm:1} follows easily by integrating by parts
\eqref{eq:standardcomp} in intervals of the form $[a,t_1]$ and
$[t_0,b]$.

Conversely, if conditions~\eqref{itm:1}, \eqref{itm:2} and
\eqref{itm:3} are satisfied, equality \eqref{eq:standardcomp}
follows easily, which implies that $\mathrm d\mathcal
L_\gamma(v)=0$ for all $v\in
T_\gamma\Omega_{PQ}\big([a,b],M;U\big)$ with small support. Since
such $v$'s span $T_\gamma\Omega_{PQ}\big([a,b],M;U\big)$, it
follows that $\gamma$ is a critical point of $\mathcal L$.
\end{proof}

We now pass to the study of the Hamiltonian formalism, and we
consider the {\em canonical symplectic form\/} $\omega$ on $TM^*$,
given by $\omega=-\mathrm d\vartheta$, where the {\em canonical
$1$-form\/} $\vartheta$ on $TM^*$ is defined by
$\vartheta_p(\zeta)=p\big(\mathrm d\pi_p(\zeta)\big)$, for all
$p\in TM^*$, $\zeta\in T_pTM^*$. If $q=(q_1,\ldots,q_n)$ is a
chart on $M$ and $(q,p)=(q_1,\ldots,q_n,p_1,\ldots,p_n)$ is the
corresponding chart on $TM^*$, the forms $\vartheta$ and $\omega$
are given by:
\begin{equation}\label{eq:varthetaomega}
\vartheta=\sum_{i=1}^np_i\,\mathrm dq_i,\quad\omega=\sum_{i=1}^n
\mathrm dq_i\wedge\mathrm dp_i.
\end{equation}
Given a Hamiltonian $H$ on $M$,
we define its {\em Hamiltonian vector field\/} $\vec H$ to be the
unique time-dependent vector field on $TM^*$ satisfying:
\[\omega(\vec H,\cdot)=\mathrm dH_t,\]
where $H_t=H(t,\cdot)$.

We say that a curve $\gamma:[a,b]\to M$ is a {\em solution of the
Hamiltonian $H$\/} if there exists a $C^1$-curve $\Gamma:[a,b]\to
TM^*$ with $\pi\circ\Gamma=\gamma$ and such that:
\begin{equation}\label{eq:Gamma}
\frac{\mathrm d}{\mathrm dt}\,\Gamma(t)=\vec
H\big(t,\Gamma(t)\big)
\end{equation} for all $t$. In this
case, we say that $\Gamma$ is a {\em Hamiltonian lift\/} of $\gamma$.
In coordinates $(q,p)$, equation \eqref{eq:Gamma} is written as:
\begin{equation}\label{eq:HJ}
\left\{\begin{aligned}
\frac{\mathrm dq}{\mathrm dt}=\;&\frac{\partial H}{\partial
p}\big(t,q(t),p(t)\big),\\[3pt]
\frac{\mathrm dp}{\mathrm dt}=\;&{-\frac{\partial H}{\partial
q}\big(t,q(t),p(t)\big)}.
\end{aligned}\right.
\end{equation}
These are called the {\em Hamilton equations\/} of $H$; observe
that the first equation in \eqref{eq:HJ} can be written
intrinsically as:
\begin{equation}\label{eq:intrinsic}
\dot\gamma(t)=\mathbb FH^{(2)}\big(t,\Gamma(t)\big).
\end{equation}
\begin{teo}\label{thm:LagrHamil}
Let $L$ be a hyper-regular Lagrangian on $M$ and let $H=L^*$ be
the corresponding hyper-regular Hamiltonian. Let $P$ and $Q$ be
smooth submanifolds of $M$; a curve
$\gamma\in\Omega_{PQ}\big([a,b],M;U\big)$ is a critical point of
$\mathcal L$ if and only if $\gamma$ is a solution of the
Hamiltonian $H$ which admits a Hamiltonian lift $\Gamma$ such that
\begin{equation}\label{eq:Gammaanula}
\Gamma(a)\vert_{T_{\gamma(a)}P}=0,\quad
\Gamma(b)\vert_{T_{\gamma(b)}Q}=0.
\end{equation}
\end{teo}
\begin{proof}
Let $\gamma\in\Omega_{PQ}\big([a,b],M;U\big)$ be a critical point
of $\mathcal L$; set $\Gamma(t)=\mathbb
FL^{(2)}\big(t,\dot\gamma(t)\big)$. Since $\mathbb FH$ and $\mathbb
FL$ are mutually inverse, equation \eqref{eq:intrinsic} follows.
Moreover, by Proposition~\ref{thm:EL}, $\Gamma$ is of class $C^1$
and \eqref{eq:Gammaanula} holds. We now prove that the second
Hamilton equation holds, using a chart $(q,p)$ of $TM^*$. To this
aim, we differentiate with respect to $q$ the equality:
\[H\Big(t,q,\frac{\partial L}{\partial \dot q}(t,q,\dot q)\Big)=\frac{\partial L}{\partial\dot q}(t,q,\dot q)\,\dot q-L(t,q,\dot q),\]
obtaining:
\begin{equation}\label{eq:deldelq}
\frac{\partial H}{\partial q}(t,q,p)+\frac{\partial H}{\partial p}(t,q,p)\,
\frac{\partial^2 L}{\partial q\,\partial\dot q}(t,q,\dot q)=
\frac{\partial^2 L}{\partial q\,\partial\dot q}(t,q,\dot q)\,\dot q-
\frac{\partial L}{\partial q}(t,q,\dot q),
\end{equation}
where $p=\frac{\partial L}{\partial\dot q}(t,q,\dot q)$. Using
that $\mathbb FH$ and $\mathbb FL$ are mutually inverse, we get
$\frac{\partial H}{\partial p}(t,q,p)=\dot q$; it follows from
\eqref{eq:deldelq} that:
\begin{equation}\label{eq:precisa}
\frac{\partial H}{\partial q}(t,q,p)=-\frac{\partial L}{\partial
q}(t,q,\dot q).
\end{equation}
The second Hamilton equation now follows from formula
\eqref{eq:precisa} and from the Euler--Lagrange equation
\eqref{eq:EL}.

Conversely, suppose that $\gamma$ is a solution of the Hamiltonian
$H$ which admits a Hamiltonian lift $\Gamma$ satisfying
\eqref{eq:Gammaanula}. Since $\mathbb FH$ and $\mathbb FL$ are
mutually inverse, from \eqref{eq:intrinsic} it follows that
$\Gamma(t)=\mathbb FL^{(2)}\big(t,\dot\gamma(t)\big)$. Finally,
equality \eqref{eq:precisa} and the second Hamilton equation imply
the Euler--Lagrange equation \eqref{eq:EL}, and the conclusion
follows from Proposition~\ref{thm:EL}.
\end{proof}
\end{section}

\begin{section}{The Space of Horizontal Curves and its Differentiable Structure}
\label{sec:horizontal}

In this section we recall some results concerning the manifold
structure of the set of horizontal curves connecting two fixed
submanifolds of a given manifold. Most of the material presented here is well-known in
the context of sub-Riemannian geometry (see \cite{Bis, BH, LS,
Mo1, Mo2, Mo3}). Detailed proofs can be found in \cite{PTJGP}.
Actually, some minor adaptations of the proofs of \cite{PTJGP}
have to be made due to the fact that \cite{PTJGP} deals with
curves of Sobolev class $H^1$ while we have to deal here\footnote{%
This is due to the fact that the {\em sub-Riemannian energy
functional\/} studied in \cite{PTJGP} is smooth on the space of
$H^1$ curves while the action functional of an arbitrary
Lagrangian is not in general even well-defined on such space.}
with curves of class $C^1$.

Throughout the section we consider fixed an $n$-dimensional differentiable
manifold $M$ and a smooth distribution $\D\subset TM$ on $M$ of rank $k\le
n$. By a {\em horizontal curve\/} we mean a curve $\gamma:[a,b]\to M$ of class
$C^1$ with $\gamma'(t)\in\D$ for all $t\in[a,b]$. Given smooth embedded
submanifolds $P,Q\subset M$ we consider the following spaces:
\begin{align*}
\Omega\big([a,b],M\big)&=\big\{\gamma:[a,b]\to M:\text{$\gamma$ is of class
$C^1$}\big\};\\
\Omega_P\big([a,b],M\big)&=\big\{\gamma\in\Omega\big([a,b],M\big):\gamma(a)\in
P\big\};\\
\Omega_{PQ}\big([a,b],M\big)&=\big\{\gamma\in\Omega\big([a,b],M\big):\gamma(a)\in
P,\ \gamma(b)\in Q\big\};\\
\Omega\big([a,b],M,\D\big)&=\big\{\gamma\in\Omega\big([a,b],M\big):\text{$\gamma$
is horizontal}\big\};\\
\Omega_P\big([a,b],M,\D\big)&=\Omega_P\big([a,b],M\big)\cap\Omega\big([a,b],M,\D\big);\\
\Omega_{PQ}\big([a,b],M,\D\big)&=\Omega_{PQ}\big([a,b],M\big)\cap\Omega\big([a,b],M,D\big).
\end{align*}
It is well-known that $\Omega\big([a,b],M\big)$ has a natural
structure of a Banach manifold (see for instance \cite{Palais,
PT1}) and that $\Omega_P\big([a,b],M\big)$ and
$\Omega_{PQ}\big([a,b],M\big)$ are embedded Banach submanifolds of
$\Omega\big([a,b],M\big)$. Also $\Omega_P\big([a,b],M,\D\big)$ is
an embedded Banach submanifold of $\Omega\big([a,b],M\big)$. The
proof of this fact is obtained by using a suitable atlas for
$\Omega\big([a,b],M\big)$ whose construction is described below.

If $\xi$ is a vector bundle over $M$ then a {\em time-dependent
referential\/} of $\xi$ over an open subset $A\subset\R\times M$
is a family $(X_i)_{i=1}^k$ of smooth maps $X_i:A\to\xi$ such that
$\big(X_i(t,m)\big)_{i=1}^k$ is a basis of the fiber $\xi_m$ for
all $(t,m)\in A$. Given a time-dependent referential
$(X_i)_{i=1}^n$ of the tangent bundle $TM$ over an open subset
$A\subset\R\times M$, we define a map:
\[\mathcal B:\Omega\big([a,b],M;\hat A\big)\longrightarrow
C^0\big([a,b],\R^n\big),\] by $\mathcal
B(\gamma)=(h_1,\ldots,h_n)$, where:
\[\gamma'(t)=\sum_{i=1}^nh_i(t)X_i\big(t,\gamma(t)\big),\]
for all $t\in[a,b]$ and: \begin{align} \label{eq:Ahat}\hat
A&=\big\{(t,v)\in\R\times TM:\big(t,\pi(v)\big)\in A\big\},\\
\label{eq:OmegaAhat}\Omega\big([a,b],M;\hat
A\big)&=\big\{\gamma\in\Omega\big([a,b],M\big):\big(t,\gamma'(t)\big)\in\hat
A,\ \text{for all}\ t\in[a,b]\big\}.
\end{align}

\begin{lem}\label{thm:Bismustchart}
If $\phi:U\subset M\to\widetilde U\subset\R^n$ is a local chart on
$M$ and $\mathcal B$ is defined as above then the map:
\begin{equation}\label{eq:cartadeBismut}
\big\{\gamma\in\Omega\big([a,b],M;\hat A\big):\gamma(a)\in
U\big\}\ni\gamma\longmapsto\big(\phi(\gamma(a)),\mathcal
B(\gamma)\big)\in\R^n\times C^0\big([a,b],\R^n\big),
\end{equation}
is a local chart on the Banach manifold $\Omega\big([a,b],M\big)$.
\end{lem}
\begin{proof}
It is a simple application of the Inverse Function Theorem on
Banach manifolds (see \cite[Corollary~4.2]{PTJGP} for details on a
similar construction).
\end{proof}

The proposition below implies that the local charts defined on
Lemma~\ref{thm:Bismustchart} form an atlas for
$\Omega\big([a,b],M\big)$.
\begin{prop}\label{thm:timedeprefbund}
Let $\xi$ be a vector bundle over a differentiable manifold $M$.
Given a continuous curve $\gamma:[a,b]\to M$, there exists a
time-dependent referential $(X_i)_{i=1}^k$ of $\xi$ whose domain
$A$ is an open neighborhood of the graph of $\gamma$ in $\R\times
M$, i.e., $\big(t,\gamma(t)\big)\in A$ for all $t\in[a,b]$.
\end{prop}
\begin{proof}
See \cite[Lemma~2.3]{PTJGP}.
\end{proof}

Using the atlas constructed above we can prove easily that
$\Omega_P\big([a,b],M,\D\big)$ is a submanifold of
$\Omega\big([a,b],M\big)$.
\begin{prop}\label{thm:horizontsubman}
$\Omega_P\big([a,b],M,\D\big)$ is an embedded Banach submanifold
of\/ $\Omega\big([a,b],M\big)$.
\end{prop}
\begin{proof}
Applying Proposition~\ref{thm:timedeprefbund} to the vector bundle
$\D$ and to a complementary vector bundle of $\D$ in $TM$ we
obtain a time-dependent referential $(X_i)_{i=1}^n$ of $TM$ such
that $(X_i)_{i=1}^k$ is a time-dependent referential for $\D$;
moreover, we may choose $(X_i)_{i=1}^n$ so that its domain
$A\subset\R\times M$ contains the graph of any prescribed
continuous curve in $M$. If $\phi$ is a local chart of $M$ which
sends $P$ to an open subset of
$\R^r\cong\R^r\times\{0\}\subset\R^n$ then the corresponding chart
\eqref{eq:cartadeBismut} on $\Omega\big([a,b],M\big)$ sends
$\Omega_P\big([a,b],M,\D\big)$ to an open subset of $\R^r\times
C^0\big([a,b],\R^k\big)$.
\end{proof}

Given Banach manifolds $\mathcal M$, $\mathcal N$, recall that a
map $f:\mathcal M\to\mathcal N$ of class $C^1$ is said to be a
{\em submersion\/} at a point $x\in\mathcal M$ if the differential
$\dd f_x:T_x\mathcal M\to T_{f(x)}\mathcal N$ is surjective and
its Kernel $\Ker(\dd f_x)$ is complemented in $T_x\mathcal M$,
i.e., it admits a closed complementary subspace in $T_x\mathcal
M$. When $f$ is a submersion at $x$, then the intersection of
$f^{-1}\big(f(x)\big)$ with some open neighborhood of $x$ in
$\mathcal M$ is a Banach submanifold of $\mathcal M$ whose tangent
space at $x$ is $\Ker(\dd f_x)$. More generally, if $\mathcal
P\subset\mathcal N$ is a Banach submanifold of $\mathcal N$ and
$x\in f^{-1}(\mathcal P)$ then we say that $f$ is {\em
transverse\/} to $\mathcal P$ at $x$ if the composition of $\dd
f_x$ with the quotient map $T_{f(x)}\mathcal N\to T_{f(x)}\mathcal
N/T_{f(x)}\mathcal P$ is surjective and has complemented kernel in
$T_x\mathcal M$; equivalently, $f$ is transverse to $\mathcal P$
at $x$ if $\Img(\dd f_x)+T_{f(x)}\mathcal P=T_{f(x)}\mathcal N$
and $\dd f_x^{-1}(T_{f(x)}\mathcal P)$ is complemented in
$T_x\mathcal M$. If $f$ is transverse to $\mathcal P$ at $x$ then
the intersection of $f^{-1}(\mathcal P)$ with some open
neighborhood of $x$ in $\mathcal M$ is a Banach submanifold of
$\mathcal M$ whose tangent space at $x$ is $\dd
f_x^{-1}(T_{f(x)}\mathcal P)$.

\begin{defin}\label{thm:singcurve}
A curve $\gamma\in\Omega_{PQ}\big([a,b],M,\D\big)$ is called {\em
regular\/} in $\Omega_{PQ}\big([a,b],M,\D\big)$ if the {\em
endpoint map}:
\begin{equation}\label{eq:endpointmap}
\Omega_P\big([a,b],M,\D\big)\ni\mu\longmapsto\mu(b)\in M
\end{equation}
is transverse to $Q$ at the point $\gamma$. When $\gamma$
is not regular in $\Omega_{PQ}\big([a,b],M,\D\big)$, we say that
$\gamma$ is {\em singular\/} in $\Omega_{PQ}\big([a,b],M,\D\big)$.
\end{defin}
Since $M$ is finite-dimensional, a curve $\gamma$ is regular in
$\Omega_{PQ}\big([a,b],M,\D\big)$ if and only if the image of the
differential of \eqref{eq:endpointmap} at $\gamma$ plus
$T_{\gamma(b)}Q$ equals $T_{\gamma(b)}M$.

Below we described an explicit method for computing the image of
the differential of the endpoint map.

\begin{defin}\label{thm:defcharacteristic}
Denote by $\D^\oo\subset TM^*$ the annihilator of $\D$. A curve
$\eta:[a,b]\to TM^*$ of class $C^1$ is called a {\em
characteristic\/} for $\D$ if $\eta\big([a,b]\big)\subset\D^\oo$
and $\eta'(t)\in T_{\eta(t)}\D^\oo$ belongs to the kernel of the
restriction of $\omega_{\eta(t)}$ to $T_{\eta(t)}\D^\oo$ (recall \eqref{eq:varthetaomega}).
\end{defin}

\begin{prop}\label{thm:calcimageend}
The annihilator of the image of the differential of
\eqref{eq:endpointmap} at a curve $\gamma$ is the subspace of\/
$T_{\gamma(b)}M^*$ given by:
\[\big\{\eta(b):\text{$\eta$ is a characteristic of\/ $\D$},\
\pi\circ\eta=\gamma,\ \eta(a)\vert_{T_{\gamma(a)}P}=0\big\}.\]
\end{prop}
\begin{proof}
The proof is a minor adaptation of the proof of
\cite[Theorem~4.9]{PTJGP} where we consider the case that $P$ is a
point and we use $H^1$ curves instead of $C^1$ curves.
\end{proof}

\begin{cor}\label{thm:outrosing}
A curve $\gamma\in\Omega_{PQ}\big([a,b],M,\D\big)$ is singular in
$\Omega_{PQ}\big([a,b],M,\D\big)$ if and only if there
exists a non zero characteristic $\eta:[a,b]\to TM^*$ of\/ $\D$
with $\pi\circ\eta=\gamma$ and $\eta(a)\vert_{T_{\gamma(a)}P}=0$,
$\eta(b)\vert_{T_{\gamma(b)}Q}=0$.
\end{cor}
\begin{proof}
Follows from Proposition~\ref{thm:calcimageend} observing that a
characteristic $\eta:[a,b]\to TM^*$ that vanishes at some
$t_0\in[a,b]$ is identically zero (see \cite[Lemma~4.8]{PTJGP}).
\end{proof}

\begin{cor}\label{thm:nosingular}
If either $T_{\gamma(a)}P+\D_{\gamma(a)}=T_{\gamma(a)}M$ or
$T_{\gamma(b)}Q+\D_{\gamma(b)}=T_{\gamma(b)}M$ then $\gamma$ is
regular in $\Omega_{PQ}\big([a,b],M,\D\big)$.\qed
\end{cor}

\begin{rem}\label{thm:stronglybracket}
If the distribution $\D$ satisfies a strong non integrability
condition (for instance, if $\D$ is a {\em contact distribution})
then the restriction of the symplectic form $\omega$ to the
annihilator $\D^\oo$ of $\D$ is nondegenerate outside the zero
section and therefore all non zero characteristic curves of $\D$
are constant. In particular, every non constant curve in
$\Omega_{PQ}\big([a,b],M,\D\big)$ is regular.
\end{rem}

So far we have looked at $\Omega_{PQ}\big([a,b],M,\D\big)$ as the
set of curves $\gamma$ in the Banach manifold
$\Omega_P\big([a,b],M,\D\big)$ satisfying the constraint
$\gamma(b)\in Q$. We could also think of
$\Omega_{PQ}\big([a,b],M,\D\big)$ as the set of curves in the
Banach manifold $\Omega_{PQ}\big([a,b],M\big)$ satisfying the
constraint $\Img(\gamma')\subset\D$. Actually, the latter point of
view will be needed in the proof of our main theorem in
Subsection~\ref{sub:proof}. Our goal now is to show that both
constraints have the same singularities. This fact was shown in
\cite{PTJGP} in the context of curves of class $H^1$. However, in
the case of curves of class $C^1$ the problem is a little harder
due to the fact that not every closed subspace of a Banach space
is complemented. We have thus decided to give all the details of
the proof.

The lemma below is a general principle that says that if a set is
defined by two constraints then the singularities of the first in
the space defined by the second constraint equals the
singularities of the second in the space defined by the first.

\begin{lem}\label{thm:doisvinculos}
Let $\mathcal M$, $\mathcal N_1$, $\mathcal N_2$ be Banach
manifolds and $\mathcal P_1\subset\mathcal N_1$, $\mathcal
P_2\subset\mathcal N_2$ be Banach submanifolds. Assume that we are
given maps $f_i:\mathcal M\to\mathcal N_i$, $i=1,2$, of class
$C^1$ and a point $x\in f_1^{-1}(\mathcal P_1)\cap
f_2^{-1}(\mathcal P_2)$ such that $f_i$ is transverse to $\mathcal
P_i$ at $x$, $i=1,2$. Then the restriction
$f_1\vert_{f_2^{-1}(\mathcal P_2)}$ is transverse to $\mathcal
P_1$ at $x$ if and only if the restriction
$f_2\vert_{f_1^{-1}(\mathcal P_1)}$ is transverse to $\mathcal
P_2$ at $x$.
\end{lem}
\begin{proof}
Consider the Banach spaces $X=T_x\mathcal M$,
$Y_i=T_{f_i(x)}\mathcal N_i/T_{f_i(x)}\mathcal P_i$, $i=1,2$, and
the continuous linear maps $L_i:X\to Y_i$, $i=1,2$, given by
composition of $\dd f_i(x)$ with the quotient map
$T_{f_i(x)}\mathcal N_i\to T_{f_i(x)}\mathcal
N_i/T_{f_i(x)}\mathcal P_i$. We know that both $L_1$ and $L_2$ are
surjective and have complemented kernel. We have to show that
$L_1\vert_{\Ker(L_2)}$ is surjective with complemented kernel if
and only if $L_2\vert_{\Ker(L_1)}$ is surjective with complemented
kernel. To this aim, observe first that $L_1\vert_{\Ker(L_2)}$ is
surjective if and only if $\Ker(L_1)+\Ker(L_2)=X$ and the latter
condition is symmetric in $L_1$ and $L_2$. Finally, to complete
the proof we show that, given $i=1,2$, then
$\Ker(L_1)\cap\Ker(L_2)$ is complemented in $\Ker(L_i)$ if and
only if it is complemented in $X$. If $\Ker(L_1)\cap\Ker(L_2)$ is
complemented in $X$ then by intersecting a closed complement of
$\Ker(L_1)\cap\Ker(L_2)$ in $X$ with $\Ker(L_i)$ we obtain a
closed complement of $\Ker(L_1)\cap\Ker(L_2)$ in $\Ker(L_i)$.
Conversely, if $Z$ is a closed complement of
$\Ker(L_1)\cap\Ker(L_2)$ in $\Ker(L_i)$ and $Z'$ is a closed
complement of $\Ker(L_i)$ in $X$ then $Z\oplus Z'$ is a closed
complement of $\Ker(L_1)\cap\Ker(L_2)$ in $X$ because
$X=\Ker(L_i)\oplus Z'$ has the product topology of $\Ker(L_i)$ and
$Z'$.
\end{proof}

We can now prove the following:
\begin{prop}\label{thm:anotherconstraint}
Let $(\theta_i)_{i=1}^{n-k}$ be a time-dependent referential of
$\D^\oo$ defined over an open subset $A\subset\R\times M$; set
$\theta=(\theta_1,\ldots,\theta_{n-k})$, so that
$\theta_{(t,m)}:T_mM\to\R^{n-k}$ is a surjective linear map with
$\Ker(\theta_{(t,m)})=\D_m$ for all $(t,m)\in A$. Consider the
map:
\[\Theta:\Omega_{PQ}\big([a,b],M;\hat A\big)\longrightarrow
C^0\big([a,b],\R^{n-k}\big)\] defined by:
\[\Theta(\gamma)(t)=\theta\big(\gamma'(t)\big),\quad t\in[a,b].\]
Then $\gamma$ is regular in $\Omega_{PQ}\big([a,b],M,\D\big)$
(in the sense of Definition~\ref{thm:singcurve}) if and only if\/
$\Theta$ is a submersion at $\gamma$.
\end{prop}
\begin{proof}
Let $\overline\Theta$ denote the extension of $\Theta$ to
\[\Omega_P\big([a,b],M;\hat
A\big)=\Omega_P\big([a,b],M\big)\cap\Omega\big([a,b],M;\hat
A\big)\]
which is again defined by
$\overline\Theta(\gamma)(t)=\theta\big(\gamma'(t)\big)$. The
conclusion will follow by applying Lemma~\ref{thm:doisvinculos}
with $\mathcal M=\Omega_P\big([a,b],M;\hat A\big)$, $\mathcal
N_1=C^0\big([a,b],\R^{n-k}\big)$, $\mathcal P_1=\{0\}$, $\mathcal
N_2=M$, $\mathcal P_2=Q$, $f_1=\overline\Theta$ and $f_2:\mathcal
M\to\mathcal N_2$ equal to the endpoint map $\mu\mapsto\mu(b)$.
Since $f_2$ is obviously a submersion, we only need to show that
$f_1=\overline\Theta$ is a submersion. Choose a distribution
$\D'\subset TM$ with $TM=\D\oplus\D'$ and let $(X_i)_{i=1}^{n-k}$
be the time-dependent referential of $\D'$ over $A$ which is dual
to $(\theta_i)_{i=1}^{n-k}$, i.e., $\theta_i(X_j)=1$ for $i=j$ and
$\theta_i(X_j)=0$ for $i\ne j$. Choose a time-dependent
referential $(X_i)_{i=n-k+1}^n$ of $\D$ over an open neighborhood
of the graph of $\gamma$. The coordinate representation of
$\Theta$ in the chart \eqref{eq:cartadeBismut} corresponding to
$(X_i)_{i=1}^n$ is the natural projection of $\R^n\oplus
C^0\big([a,b],\R^n\big)$ onto $C^0\big([a,b],\R^{n-k}\big)$. This
shows that $\overline\Theta$ is a submersion and concludes the
proof.
\end{proof}

\end{section}

\begin{section}{Lagrangians with linear constraints and degenerate Hamiltonians}
\label{sec:degenerate} Let $M$ be an $n$-dimensional manifold and
$\D\subset TM$ be a smooth distribution of rank $k$. We consider
$\D$ as a vector bundle over $M$ with projection $\pi:\D\to M$. We
apply the theory of Subsection~\ref{sub:legendre} to the vector
bundle $\xi=\R\times\D$ over the manifold $\R\times M$, with
projection $\mathrm{Id}\times\pi$. The fiber $\xi_{(t,m)}$ is
given by $\{t\}\times\D_m$.

Let $L:U\subset \xi\to\R$ be a map of class $C^2$ defined in an
open set $U\subset\xi$; we assume that $L$ is hyper-regular in the
sense of Definition~\ref{thm:defFZ}, so that (by the Inverse
Function Theorem) the fiber derivative $\mathbb FL:U\to V$ is a
$C^1$ diffeomorphism onto an open subset $V\subset\xi^*$. Let
$H_0=L^*$ be the Legendre transform of $L$. Then $H_0:V\to\R$ is a
map of class $C^1$ whose restriction to each fiber of $\xi^*$ is
of class $C^2$; moreover, the fiber derivative $\mathbb FH_0:V\to
U$ is the inverse of $\mathbb FL$ (see
Proposition~\ref{thm:bundles}).

For every $p\in TM^*$ we denote by $p\vert_{\D}$ the restriction of $p\in
T_mM^*$ to $\D_m$. Observe that the {\em restriction map\/} $TM^*\ni p\to
p\vert_{\D}\in\D^*$ is the transpose of the vector bundle inclusion $\D\to
TM$. By composing $H_0$ with the restriction map $TM^*\to\D^*$ we obtain a map
$H:\widetilde V\to\R$ given by:
\begin{equation}\label{eq:asin}
H(t,p)=H_0\big(t,p\vert_{\D}\big),\quad(t,p)\in\widetilde V,
\end{equation}
where:
\[\widetilde V=\big\{(t,p)\in\R\times TM^*:(t,p\vert_{\D})\in V\big\}.\]
Observe that $H$ is a Hamiltonian on $M$ (see Definition~\ref{thm:defLagrHam})
of class $C^1$ defined in the
open set $\widetilde V\subset\R\times TM^*$. We will call $L$ a {\em
constrained Lagrangian\/} on $M$, and $H$ the corresponding {\em degenerate
Hamiltonian\/} (observe indeed that $H$ cannot be regular unless $\D=TM$).

Given any two submanifolds $P$ and $Q$ of $M$ then a constrained
Lagrangian $L$ on $M$ defines an action functional $\mathcal L$ on
$\Omega_{PQ}\big([a,b],M,\D;U\big)$ by formula \eqref{eq:action}.
Our goal is to determine the critical points of $\mathcal L$.

The following is the main result of the paper and its proof is given in
Subsection~\ref{sub:proof}:
\begin{teo}\label{thm:degenerate}
Let $M$ be an $n$-dimensional manifold, $\D\subset TM$ be a smooth
distribution of rank $k$ and $L:U\subset\R\times\D\to\R$ be a
hyper-regular constrained Lagrangian of class $C^2$. Let $H_0=L^*$
be the Legendre transform of $L$ and let $H$ be the corresponding
degenerate Hamiltonian as in \eqref{eq:asin}.

Fix two submanifolds $P$ and $Q$ of $M$ and let $\mathcal L$ be
the action functional of $L$ defined in the space
$\Omega_{PQ}\big([a,b],M,\D;U\big)=\Omega_{PQ}\big([a,b],M,\D\big)\cap\Omega_{PQ}\big([a,b],M;U\big)$,
given by \eqref{eq:action}. Let
$\gamma\in\Omega_{PQ}\big([a,b],M,\D;U\big)$ be a regular curve.
Then, $\gamma$ is a critical point of $\mathcal L$ if and only if
it is a solution of $H$ that admits a Hamiltonian lift
$\Gamma:[a,b]\to TM^*$ with $\Gamma(a)\vert_{T_{\gamma(a)}P}=0$
and $\Gamma(b)\vert_{T_{\gamma(b)}Q}=0$.
\end{teo}

The classical example of a constrained hyper-regular Lagrangian
function $L$ is given by:
\begin{equation}\label{eq:classicalL}
L(t,v)=\tfrac12\,g(v,v)-V\big(\pi(v)\big),
\end{equation}
where $g$ is a {\em sub-Riemannian metric\/} on $(M,\D)$ (i.e., a smooth Riemannian structure on the vector bundle $\D$) and $V:M\to\R$ is a map of class
$C^2$. The fiber derivative $\mathbb FL$ of \eqref{eq:classicalL} is given by:
\[\mathbb FL(t,v)=g(v,\cdot)\in\D^*,\]
so that $L$ is indeed hyper-regular. Recalling \eqref{eq:EZfibrado}, we compute as follows:
\begin{gather*}
E_L(t,v)=\tfrac12\,g(v,v)+V\big(\pi(v)\big),\quad v\in\D,\\
H_0(t,\rho)=\tfrac12\,g^{-1}(\rho,\rho)+V\big(\pi(\rho)\big),\quad
\rho\in\D^*,
\end{gather*}
where $g^{-1}$ denotes the induced Riemannian structure on the dual bundle
$\D^*$. The degenerate Hamiltonian $H$ corresponding to \eqref{eq:classicalL}
is thus given by:
\[H(t,p)=\tfrac12\,g^{-1}(p\vert_\D,p\vert_\D)+V\big(\pi(p)\big),\quad p\in
TM^*.\]
Theorem~\ref{thm:degenerate} implies that the critical points of the
action functional $\mathcal L$ corresponding to \eqref{eq:classicalL} on the space
$\Omega_{PQ}\big([a,b],M,\D\big)$ are the
solutions of $H$ that admit a Hamiltonian lift $\Gamma:[a,b]\to TM^*$
satisfying the boundary conditions $\Gamma(a)\vert_{T_{\gamma(a)}P}=0$ and
$\Gamma(b)\vert_{T_{\gamma(b)}Q}=0$. Observe that (in the case when
$P$ and $Q$ are points) we obtain the equations for the trajectories of the
Vakonomic mechanics given in \cite{KO}; when $V=0$ we obtain the equations for
the normal geodesics of the sub-Riemannian manifold $(M,\D,g)$ (see
\cite{LS}).

\begin{rem}\label{thm:minimoanormal}
We emphasize that, in general, a minimum of the action functional
$\mathcal L$ may not be a regular curve in
$\Omega_{PQ}\big([a,b],M,\D\big)$, and in this situation it may
not satisfy the Hamilton equations of $H$. Examples of this
phenomenon are given in \cite{LS, Mo1} in the sub-Riemannian case
$L(t,v)=\frac12g(v,v)$. Hence, one can only conclude that a
minimum of $\mathcal L$ is either a solution of the Hamilton
equations or the projection of a non null characteristic of $\D$.
\end{rem}

\mysubsection{Generalized functions calculus}\label{sub:calculus}
For the proof of Theorem~\ref{thm:degenerate} we will occasionally have to consider
derivatives of functions that are in principle only continuous\footnote{%
This situation already occurred in the proof of Proposition~\ref{thm:EL}. In that
case the difficulty could also be circumvented by a simpler technique.}. These
derivatives should be understood in the sense of
Schwartz's generalized functions calculus.
However, the usual definition of the generalized functions
space as the dual of the space of smooth compactly supported
maps only allows products of generalized functions by smooth
maps. To overcome this difficulty, we introduce a calculus
for generalized functions of {\em stronger\/} regularity,
that are elements of the dual of a space of functions with
{\em weaker\/} regularity.

Let $V$ be a real finite dimensional vector space. For $k\ge0$, we define
$C^k_0\big([a,b],V\big)$ to be the Banach space of $V$-valued $C^k$ maps on $[a,b]$
whose first $k$ derivatives vanish at $a$ and at $b$; we endow it with the standard
$C^k$-norm. We denote by $D^k\big([a,b],V\big)$ the dual Banach space of
$C^k_0\big([a,b],V^*\big)$ (dual spaces will {\em always\/} be meant in the
{\em topological sense}). Denoting by $L^p\big([a,b],V\big)$ the Banach space
of $V$-valued
measurable functions on $[a,b]$ whose $p$-th power is Lebesgue integrable, we have an inclusion:
\begin{equation}\label{eq:inclL1Dk}
L^1\big([a,b],V\big)\hookrightarrow D^k\big([a,b],V\big)
\end{equation}
defined by \[\langle f,\alpha\rangle=\int_a^b\alpha(t)\,f(t)\;\mathrm dt,\quad
f\in L^1\big([a,b],V\big),\ \alpha\in C^k_0\big([a,b],V^*\big);\]
in the formula above we have denoted by $\langle f,\alpha\rangle$ the
evaluation at $\alpha$ of the linear functional which is the image of $f$ by
\eqref{eq:inclL1Dk}. In what follows we will always identify a function $f\in
L^1\big([a,b],V\big)$ with its image by \eqref{eq:inclL1Dk}; moreover, the
evaluation of $f\in D^k\big([a,b],V\big)$ at $\alpha\in
C^k_0\big([a,b],V^*\big)$ will always be denoted by $\langle f,\alpha\rangle$.

Observe that we have inclusions $D^k\hookrightarrow D^{k+1}$ defined by
restriction of the functionals, i.e., $D^k\hookrightarrow D^{k+1}$ is simply
the transpose of the inclusion of $C^{k+1}_0\big([a,b],V^*\big)$ in
$C^k_0\big([a,b],V^*\big)$.

We summarize the observations above by the following
diagram:
\[\cdots\hookrightarrow C^1 \hookrightarrow C^0\hookrightarrow
L^1\hookrightarrow D^0 \hookrightarrow D^1\hookrightarrow\cdots\]
An element $f$ of any space $D^k\big([a,b],V\big)$ is called a {\em
generalized function}. In what follows, we will occasionally write simply
$C^k$, $C^k_0$, $D^k$, $L^p$ instead of $C^k\big([a,b],V\big)$,
$C^k_0\big([a,b],V\big)$, $D^k\big([a,b],V\big)$, $L^p\big([a,b],V\big)$.

In addition to the standard vector space operations in $D^k$, we define
the following:
\begin{itemize}
\item {\em derivative operation:} for $f\in D^k\big([a,b],V\big)$, we define
the {\em derivative\/} of $f$ to be the generalized function
$f'\in D^{k+1}\big([a,b],V\big)$ defined by:
\[\langle f',\alpha\rangle=-\langle f,\alpha'\rangle,\]
for all $\alpha\in C^{k+1}_0\big([a,b],V^*\big)$;
\smallskip

\item {\em product operation:} for $f\in D^k\big([a,b],V\big)$, $g\in C^k\big([a,b],W\big)$ and
a fixed bilinear map $V\times W\to U$, we define the product $fg\in
D^k\big([a,b],U\big)$ as follows. The bilinear map $V\times W\to U$
induces a bilinear map $W\times U^*\to V^*$ defined by $(w\cdot
u^*)(v)=u^*(v\cdot w)$; we set:
\[\langle fg,\alpha\rangle=\langle f,g\cdot\alpha\rangle,\]
for all $\alpha\in C^k_0\big([a,b],U^*\big)$;
\smallskip

\item {\em restriction operation:} for $f\in D^k\big([a,b],V\big)$ and $[c,d]\subset[a,b]$, we set:
\[\langle f\vert_{[c,d]},\alpha\rangle=\langle f,\overline\alpha\rangle,\]
for all $\alpha\in C^k_0\big([c,d],V^*\big)$, where $\overline\alpha\in
C^k_0\big([a,b],V^*\big)$ is the extension to zero
of $\alpha$ outside $[c,d]$.

\end{itemize}

It is easily seen that when we apply the above operations to elements of $D^k$
which correspond to functions then we obtain the standard operations on
functions. Moreover, the standard Leibnitz rule for
derivatives of products holds for generalized functions, i.e.:
\[(fg)'=f'g+fg',\]
for all $f\in D^k$ and $g\in C^{k+1}$.

In order to prove some regularity results we present the following elementary
lemmas.

\begin{lem}\label{thm:dergeneralizada}
Let $f\in D^k\big([a,b],V\big)$ be such that $f'=0$. Then $f$ is a constant function.
\end{lem}
\begin{proof}
We first consider the case $V=\R$. If $f'=0$, then $\langle f,\alpha'\rangle=0$ for
all $\alpha\in C^{k+1}_0\big([a,b],\R\big)$, hence $\langle f,\beta\rangle=0$
for all
$\beta\in C^k_0\big([a,b],\R\big)$ with $\int_a^b\beta=0$. Choose $\beta_0\in
C^k_0\big([a,b],\R\big)$ with $\int_a^b\beta_0=1$; set $c=\langle
f,\beta_0\rangle$. It is easily seen that $f\equiv c$.

For the general case, observe that for all $\lambda\in V^*$, the product
$\lambda\,f\in D^k\big([a,b],\R\big)$ has vanishing derivative, and hence it
is constant. Since $\lambda$ is arbitrary, it follows that $f$ is constant.
\end{proof}

\begin{lem}\label{thm:existenciaprimitiva}
Let $f\in D^k\big([a,b],V\big)$, $k\ge1$; there exists an element $F\in
D^{k-1}\big([a,b],V\big)$
with
$F'=f$. If $f\in D^0\big([a,b],V\big)$, there exists $F\in
L^2\big([a,b],V\big)$ with $F'=f$.
\end{lem}
\begin{proof}
Consider the map $\mathrm d:C^{k+1}_0\to C^k_0$ given by $\mathrm
d(\alpha)=\alpha'$. It is easily seen that $\mathrm d$ is injective with
closed image. It follows that the {\em transpose
map\/} $\mathrm d^*:D^k\to D^{k+1}$ is surjective; clearly, the
derivative operator for generalized functions is $-\mathrm d^*$, which proves the
first part of the thesis.

For the case $k=0$, let $H^1_0$ denote the Sobolev space of
absolutely continuous functions $\alpha:[a,b]\to V^*$ having
square integrable derivative, and such that
$\alpha(a)=\alpha(b)=0$. Again, the derivation map $\mathrm
d:H^1_0\to L^2$ is injective and has closed
image. Therefore, given $f\in D^0$, we can find $F\in
{L^2}^*\simeq L^2$ with $\mathrm d^*F=-f\vert_{H^1_0}$. It follows
that $F'=f$.
\end{proof}

\begin{cor}[Bootstrap lemma]\label{thm:bootstrap}
Let $f$ be a generalized function.

\begin{enumerate}
\item If $f'\in D^0$ then $f\in L^2$;

\item If $f'\in L^2$ then $f\in C^0$;

\item If $f'\in C^0$ then $f\in C^1$.
\end{enumerate}
\end{cor}
\begin{proof}
We prove, for example, the first item. By Lemma~\ref{thm:existenciaprimitiva},
we can find $F\in L^2$ with $F'=f'$. By Lemma~\ref{thm:dergeneralizada}, it
follows that $F-f$ is constant, hence $f\in L^2$.

The other items are proven similarly.
\end{proof}

We now give a result that shows that {\em regularity\/}
of a generalized function is a local property:
\begin{lem}\label{thm:regloc}
Let $\lambda$ be a generalized function on $[a,b]$. Suppose that
for all $t\in[a,b]$ there exists $\varepsilon>0$ such that the
restriction $\lambda\vert_{[t-\varepsilon,t+\varepsilon]\cap[a,b]}$
is of class $C^k$, $k\ge0$. Then $\lambda$ is of class $C^k$.
\end{lem}
\begin{proof}
Consider a partition $a=t_0<t_1<\ldots<t_r=b$ such that $f_i=\lambda\vert_{[t_i,t_{i+2}]}$ is
of class $C^k$ for all $i=0,\ldots,r-2$. Since the operation of restriction
for generalized functions gives the standard operation of restriction for
functions, it follows that:
\[f_i\vert_{[t_{i+1},t_{i+2}]}=\lambda\vert_{[t_{i+1},t_{i+2}]}=f_{i+1}\vert_{[t_{i+1},t_{i+2}]},\]
for $i=0,\ldots,r-3$. Hence there exists a $C^k$ map $f$ on $[a,b]$ such that
$f\vert_{[t_i,t_{i+2}]}=f_i$ for all $i=0,\ldots,r-2$. We know that $\langle
f,\alpha\rangle=\langle\lambda,\alpha\rangle$ if $\alpha$ has support
contained in some interval $\left]t_i,t_{i+2}\right[$; but such $\alpha$'s
span a dense subspace of the domain of the linear functional $\lambda$ and
therefore $\lambda=f$.
\end{proof}

Finally, we need the following result that relates the dual spaces of $C^0$
and $C^0_0$.
For $t\in[a,b]$ and $\sigma\in V$, we denote by $\delta_t^\sigma\in
C^0\big([a,b],V^*\big)^*$
the {\em Dirac's delta}, defined by:
\[\langle\delta_t^\sigma,\alpha\rangle=\alpha(t)\,\sigma,\quad\alpha\in
C^0\big([a,b],V^*\big).\]
\begin{lem}\label{thm:C0C00}
If $\lambda\in C^0\big([a,b],V^*\big)^*$ vanishes identically
on $C^0_0\big([a,b],V^*\big)$ then there exist $\sigma_a$ and $\sigma_b$ in
$V$ such that:
\begin{equation}\label{eq:deltaab}
\lambda=\delta_a^{\sigma_a}+\delta_b^{\sigma_b}.
\end{equation}
\end{lem}
\begin{proof}
If $\mathcal A$ denotes the subspace of $C^0\big([a,b],V^*\big)$ consisting of
{\em affine maps\/} $\alpha(t)=Pt+Q$ then obviously:
\[C^0\big([a,b],V^*\big)=C^0_0\big([a,b],V^*\big)\oplus\mathcal A.\]
It is easy to see that we can find $\sigma_a,\sigma_b\in V$ such that both
sides of \eqref{eq:deltaab} agree on $\mathcal A$. Since both sides of
\eqref{eq:deltaab} vanish on $C^0_0\big([a,b],V^*\big)$, the conclusion follows.
\end{proof}

\mysubsection{Proof of Theorem~\ref{thm:degenerate}}
\label{sub:proof}
The proof of Theorem~\ref{thm:degenerate} is based on the method of Lagrange multipliers,
and we start with the precise statement of the result needed for our purposes.

\begin{prop}\label{thm:LagrMul}
Let $\mathcal M$ be a Banach manifold, $E$ a Banach space and
$F:\mathcal M\to\R$, $g:\mathcal M\to E$ maps of class $C^1$. Let
$p\in g^{-1}(0)$ be such that $g$ is a submersion at $p$. Then,
$p$ is a critical point for $f\vert_{g^{-1}(0)}$ if and only if
there exists $\lambda\in E^*$ such that $p$ is a critical point
for the functional $f_\lambda=f-\lambda\circ g$ in $\mathcal M$.
\end{prop}
\begin{proof}
The point $p$ is critical for $f\vert_{g^{-1}(0)}$ if and only if
$\mathrm df(p)$ vanishes on $T_pg^{-1}(0)=\mathrm{Ker}\big(\mathrm
dg(p)\big)$. The conclusion follows from elementary functional
analysis arguments.
\end{proof}
The linear functional $\lambda\in E^*$ of
Proposition~\ref{thm:LagrMul} is called the {\em Lagrange
multiplier\/} of the constrained critical point $p$; it is easily
seen that such $\lambda$ is unique. We can now prove the main
result of the section. In the argument we will need a regularity
result for a Lagrangian multiplier; such proof is postponed to
Lemma~\ref{thm:postponed}.
\begin{proof}[Proof of Theorem~\ref{thm:degenerate}]
We start by choosing an arbitrary complementary distribution $\D'$
to $\D$, i.e., a smooth distribution of rank $n-k$ in $M$ such
that $T_mM=\D_m\oplus\D'_m$ for all $m\in M$; moreover, we fix an
arbitrary smooth Riemannian structure $g$ on the vector bundle
$\D'$. Let $\pi_{\D}:TM\to\D$ and $\pi_{\D'}:TM\to\D'$ denote the
projections and define an extension $\widetilde L:\widetilde
U\subset\R\times M\to\R$ of $L$ by:
\begin{equation}\label{eq:extendedL}
\widetilde
L(t,v)=L\big(t,\pi_{\D}(v)\big)+\tfrac12\,g\big(\pi_{\D'}(v),\pi_{\D'}(v)\big),
\end{equation}
where
\[\widetilde U=\Big\{(t,v)\in\R\times TM:\big(t,\pi_{\D}(v)\big)\in U\Big\}.\]
Then $\widetilde U$ is open in $\R\times TM$ and $\widetilde L$ is
a Lagrangian on $M$ as in Definition~\ref{thm:defLagrHam}; we
denote by $\widetilde{\mathcal L}$ the corresponding action
functional in $\Omega_{PQ}\big([a,b],M;\widetilde U\big)$, defined
as in \eqref{eq:action}.

Let $\theta$, $\Theta$, $A$ and $\hat A$ be as in the statement of
Proposition~\ref{thm:anotherconstraint} (recall also
\eqref{eq:Ahat} and \eqref{eq:OmegaAhat}). Then, since $\gamma$ is
regular, the map $\Theta$ is a submersion at $\gamma$; moreover,
$\gamma$ is a critical point of $\mathcal L$ in
$\Omega_{PQ}\big([a,b],M,\D;U\big)$ if and only if it is a
critical point of $\widetilde{\mathcal L}\vert_{\Theta^{-1}(0)}$.
By the method of Lagrange multipliers
(Proposition~\ref{thm:LagrMul}), this is equivalent to the
existence of $\lambda\in C^0\big([a,b],\R^{n-k}\big)^*$ such that
$\gamma$ is a critical point of $\widetilde{\mathcal L}_\lambda=
\widetilde{\mathcal L}-\lambda\circ\Theta$ in
$\Omega_{PQ}\big([a,b],M;\hat A\cap\widetilde U\big)$.

We will prove in Lemma~\ref{thm:postponed} below that the Lagrange
multiplier $\lambda$ is of class $C^1$, i.e., that it is given by:
\begin{equation}\label{eq:lambdaC1}
\lambda(\alpha)=\int_a^b\lambda_0(t)\,\alpha(t)\;\mathrm dt,\quad\forall\,\alpha\in
C^0\big([a,b],\R^{n-k}\big),
\end{equation}
for some $C^1$ map $\lambda_0:[a,b]\to(\R^{n-k})^*$. Therefore,
$\widetilde{\mathcal L}_\lambda$ is the action functional
corresponding to the Lagrangian $\widetilde L_\lambda$ in $M$
defined by:
\begin{equation}\label{eq:oLlambda}
\widetilde L_\lambda(t,v)=\widetilde
L(t,v)-\lambda_0(t)\,\theta_{(t,m)}(v),\quad (t,v)\in\hat A\cap
\widetilde U,
\end{equation}
where $m=\pi(v)$.

We now prove that $\widetilde L$ and $\widetilde L_\lambda$ are
hyper-regular and we compute their Legendre transforms. The fiber
derivatives $\mathbb F\widetilde L$ and $\mathbb F\widetilde
L_\lambda$ are easily computed as:
\begin{gather}
\label{eq:fiberLLlambda1}\mathbb F\widetilde L(t,v)=\mathbb
FL\big(t,\pi_{\D}(v)\big)\circ\pi_{\D}+
g\big(\pi_{\D'}(v),\pi_{\D'}(\,\cdot\,)\big)\in T_mM^*,\\
\label{eq:fiberLLlambda2}\mathbb F\widetilde
L_\lambda(t,v)=\mathbb F\widetilde
L(t,v)-\lambda_0(t)\,\theta_{(t,m)}\in T_mM^*,
\end{gather}
where $m=\pi(v)$. The hyper-regularity is proven by exhibiting
explicit inverses:
\begin{equation}\label{eq:inverse}
\begin{split}
&\mathbb F\widetilde L^{-1}(t,p)=\mathbb FL^{-1}(t,p\vert_{\D})+g^{-1}(p\vert_{\D'}),\\
&\mathbb F\widetilde L_\lambda^{-1}(t,p)=\mathbb F\widetilde
L^{-1}\big(t,p+\lambda_0(t)\,\theta_{(t,m)}\big);
\end{split}
\end{equation}
by $g^{-1}$ in the above formula we mean the inverse of $g$ seen
as a linear map from $\D_m$ to $\D_m^*$.

We now compute the Legendre transforms $\widetilde H$ and
$\widetilde H_\lambda$ of $\widetilde L$ and $\widetilde
L_\lambda$ respectively. Using Definition~\ref{thm:defleg} and
equations \eqref{eq:fiberLLlambda1}, \eqref{eq:fiberLLlambda2}, we
compute easily:
\begin{equation}\label{eq:ELLl}
E_{\widetilde L_\lambda}(t,v)=E_{\widetilde
L}(t,v)=E_L\big(t,\pi_{\D}(v)\big)+\tfrac12\,
g\big(\pi_{\D'}(v),\pi_{\D'}(v)\big);
\end{equation}
and, using \eqref{eq:inverse}, we therefore obtain:
\begin{align*}
&\widetilde H(t,p)=H(t,p)+\tfrac12\,g^{-1}(p\vert_{\D'},p\vert_{\D'}),\\
&\begin{aligned}\widetilde H_\lambda(t,p)=&\,\widetilde H(t,p+\lambda_0(t)\,\theta_{(t,m)})\\
=&\,H(t,p)+\tfrac12\,g^{-1}\big((p+\lambda_0(t)\,\theta_{(t,m)})\vert_{\D'},
(p+\lambda_0(t)\,\theta_{(t,m)})\vert_{\D'}\big).
\end{aligned}
\end{align*}
We now compute the Hamilton  equations of the Hamiltonian
$\widetilde H_\lambda$ with the help of local coordinates
$(q_1,\ldots,q_n,p_1,\ldots,p_n)$ in $TM^*$ and of a local
$g$-orthonormal referential $X_1,\ldots,X_{n-k}$ of $\D'$.

We write:
\begin{equation}\label{eq:HX}
\widetilde
H_\lambda(t,p)=H(t,p)+\frac12\,\sum_{i=1}^{n-k}\big(p+\lambda_0(t)\,\theta_{(t,m)}\big)(X_i)^2,
\end{equation}
and, using \eqref{eq:HJ}, the Hamilton  equations of $\widetilde
H_\lambda$ are given by:
\begin{equation}\label{eq:HJHl}
\left\{\begin{aligned}\frac{\mathrm dq}{\mathrm dt}&\,=\frac{\partial H}{\partial p}
+\sum_{i=1}^{n-k} (p+\lambda_0\,\theta)(X_i) \,X_i,\\
\frac{\mathrm dp}{\mathrm dt}&\,=-\frac{\partial H}{\partial q}-
\sum_{i=1}^{n-k}(p+\lambda_0\,\theta)(X_i)\,\left[\lambda_0\,\frac{\partial\theta}{\partial
q}(X_i)+(p+\lambda_0\,\theta)\left(\frac{\partial X_i}{\partial q}\right)\right].\end{aligned}\right.
\end{equation}
By Theorem~\ref{thm:LagrHamil}, $\gamma$ is a critical point of
$\widetilde{\mathcal L}_\lambda$ if and only if it admits a lift
$\Gamma:[a,b]\to TM^*$ satisfying \eqref{eq:HJHl} with
$\Gamma(a)\in T_{\gamma(a)}P^o$ and $\Gamma(b)\in
T_{\gamma(b)}Q^o$.

Now, it follows easily from \eqref{eq:asin} that $\frac{\partial
H}{\partial p}$ is in $\D$; since $\gamma$ is horizontal, i.e.,
$\frac{\mathrm dq}{\mathrm dt}\in\D$, from the first equation of
\eqref{eq:HJHl} it follows that $(p+\lambda_0\,\theta)(X_i)=0$ for
all $i=1,\ldots,n-k$. Setting $(p+\lambda_0\,\theta)(X_i)=0$ in
\eqref{eq:HJHl} we obtain the Hamilton equations of $H$, which
concludes the proof.
\end{proof}
We are left with the proof of the {\em regularity\/} of the
Lagrange multiplier $\lambda$. We will use the generalized
functions calculus developed in Subsection~\ref{sub:calculus}.
\begin{lem}\label{thm:postponed}
Under the assumptions of Theorem~\ref{thm:degenerate}, using the
notations adopted in its proof, if $\gamma$ is horizontal and if, for some $\lambda\in
C^0\big([a,b],\R^{n-k}\big)^*$, it is a critical point of $\widetilde{\mathcal L}-\lambda\circ\Theta$, then
there exists a $C^1$ map $\lambda_0:[a,b]\to(\R^{n-k})^*$ such that
\eqref{eq:lambdaC1} holds.
\end{lem}
\begin{proof}
We set \[\lambda_0=\lambda\vert_{C^0_0([a,b],\R^{n-k})}\in
D^0\big([a,b],(\R^{n-k})^*\big);\] we first prove the regularity
of the generalized function $\lambda_0$. To this aim, we {\em
localize\/} the problem by considering variational vector fields
along $\gamma$ having support in the domain of a local chart
$q=(q_1,\ldots,q_n)$ in $M$.

Let $[c,d]\subset[a,b]$ be such that $\gamma\big([c,d]\big)$ is contained in the
domain of the local chart; we still denote by $\lambda_0$ the restriction of
$\lambda_0$ to $[c,d]$.

Since $\gamma$ is a critical point of $\widetilde{\mathcal
L}-\lambda\circ\Theta$, by standard computations it follows that
the following equality holds:
\begin{multline}\label{eq:standardcomp2}
\int_c^d\frac{\partial\widetilde L}{\partial q}\big(t,q(t),\dot
q(t)\big)\,v(t)+\frac{\partial\widetilde
L}{\partial \dot q}\big(t,q(t),\dot q(t)\big)\,\dot v(t)\,\mathrm dt\\
-\Big\langle\lambda_0,\left.\frac{\partial \theta}{\partial
q}\right\vert_{(t,q(t))}\!\!\!\!\!\big(v(t),\dot q(t)\big)+
\theta_{(t,q(t))}\,\dot v(t)\Big\rangle=0,
\end{multline}
for  every  vector field $v$ of class $C^1$ along $\gamma$ having
support in $\left]c,d\right[$; in the formula above we have
regarded the derivative $\frac{\partial\theta}{\partial
q}\big\vert_{(t,q(t))}$ as an $\R^{n-k}$-valued bilinear map in
$\R^n$. In terms of the local coordinates, the maps $\theta$,
$\frac{\partial\theta}{\partial q}(\cdot,\dot q)$,
$\frac{\partial\widetilde L}{\partial q}$ and
$\frac{\partial\widetilde L}{\partial \dot q}$ evaluated along
$\gamma$ will be interpreted as follows:
\begin{itemize}
\item $\displaystyle\theta\in C^1\big([c,d],\mathrm{Lin}(\R^n,\R^{n-k})\big)$;
\smallskip
\item $\displaystyle\frac{\partial\theta}{\partial q}(\cdot,\dot q)\in
C^0\big([c,d],\mathrm{Lin}(\R^n,\R^{n-k})\big)$;
\smallskip
\item $\displaystyle\frac{\partial\widetilde L}{\partial q},\  \frac{\partial\widetilde L}{\partial\dot
q}\in C^0\big([c,d],{\R^n}^*\big)$,
\end{itemize}
where $\mathrm{Lin}(\cdot,\cdot)$ denotes the space of linear maps
between two given vector spaces.

Using the definition of derivative for generalized
functions, from \eqref{eq:standardcomp2} we get:
\begin{equation}\label{eq:weak}
\Big\langle{-\Big(\frac{\partial\widetilde L}{\partial\dot
q}\Big)'}+ \frac{\partial\widetilde L}{\partial
q}-\lambda_0\,\frac{\partial\theta}{\partial q}(\cdot,\dot
q)+(\lambda_0\,\theta)',v\Big\rangle=0,
\end{equation}
for  every $C^1$ map $v:[c,d]\to\R^n$  having support in
$\left]c,d\right[$, and, by density, for every $v\in
C^1_0\big([c,d],\R^n\big)$. It follows:
\begin{equation}\label{eq:strong}
-\Big(\frac{\partial\widetilde L}{\partial\dot q}\Big)'+
\frac{\partial\widetilde L}{\partial
q}-\lambda_0\,\frac{\partial\theta}{\partial q}(\cdot,\dot q)+
\lambda_0'\,\theta+\lambda_0\,\theta'=0.
\end{equation}

Let $X_1,\ldots,X_{n-k}$ be a referential of $\D'$ along $\gamma$;
in terms of the local coordinates the $X_i$'s will be thought as
elements of $C^1\big([c,d],\R^n\big)$. Moreover, we set
\[X=(X_1,\ldots,X_{n-k})\in
C^1\big([c,d],\mathrm{Lin}(\R^{n-k},\R^n)\big),\] where the
$(n-k)$-tuple $\big(X_1(t),\ldots,X_{n-k}(t)\big)$ is identified
with the linear map that takes the $i$-th vector of the canonical
basis of $\R^{n-k}$ to $X_i(t)$.

Composing \eqref{eq:strong} with $X$, we obtain:
\begin{equation}\label{eq:strongmult}
\lambda_0'\,\theta(X)+\lambda_0\,\theta'(X)-\lambda_0\,\frac{\partial\theta}{\partial
q} (X,\dot q)+\frac{\partial \widetilde L}{\partial
q}\,X-\Big(\frac{\partial \widetilde L}{\partial \dot
q}\Big)'\,X=0.
\end{equation}
Evaluating \eqref{eq:fiberLLlambda1} at $X_i$ with $v=\gamma'$ and
using the horizontality of $\gamma$ we get:
\begin{equation}\label{eq:firstXi}
\frac{\partial\widetilde L}{\partial\dot q}\,X_i=0,\quad
i=1,\ldots,n-k;
\end{equation}
hence:
\begin{equation}\label{eq:C0!}
\Big(\frac{\partial\widetilde L}{\partial\dot q}\Big)'X=-
\frac{\partial\widetilde L}{\partial\dot q}\,X'\in
C^0\big([c,d],(\R^{n-k})^*\big).
\end{equation}
Now, considering that
$\theta(X)\in\mathrm{Lin}\big(\R^{n-k},\R^{n-k}\big)$ is
invertible, by \eqref{eq:C0!} we can write \eqref{eq:strong} in
the form:
\begin{equation}\label{eq:eddifl0}
\lambda_0'=\lambda_0\,h_1+h_2,
\end{equation}
with $h_1\in C^0\big([c,d],\mathrm{Lin}(\R^{n-k},\R^{n-k})\big)$
and $h_2\in C^0\big([c,d],(\R^{n-k})^*\big)$.

Applying three times Corollary~\ref{thm:bootstrap}, from
\eqref{eq:eddifl0} we conclude that $\lambda_0$ belongs to the
space $ C^1\big([c,d],(\R^{n-k})^*\big)$; now
Lemma~\ref{thm:regloc} implies that $\lambda_0\in
C^1\big([a,b],(\R^{n-k})^*\big)$.

By Lemma~\ref{thm:C0C00}, there exist $\sigma_a,\sigma_b\in(\R^{n-k})^*$ such that:
\begin{equation}\label{eq:ll0}
\lambda(\alpha)=\int_a^b\lambda_0\,\alpha\,\mathrm
dt+\sigma_a\,\alpha(a)+\sigma_b\,\alpha(b),\quad\alpha\in
C^0\big([a,b],\R^{n-k}\big).
\end{equation}
To conclude the proof we show that $\sigma_a=\sigma_b=0$. Let's
show for instance that $\sigma_a=0$; the proof of the equality
$\sigma_b=0$ is analogous.

Using local charts around $\gamma\big([a,d]\big)$, for $d$ close
to $a$, we consider variational vector fields $v$ of class $C^1$
supported in $\left[a,d\right[$, with $v(a)\in T_{\gamma(a)}P$.
Arguing as in the deduction of formula \eqref{eq:standardcomp2},
we get the following equality:
\begin{multline}\label{eq:standardcomp3}
\int_a^d\frac{\partial \widetilde L}{\partial q}\big(t,q(t),\dot
q(t)\big)\,v(t)+
\frac{\partial \widetilde L}{\partial \dot q}\big(t,q(t),\dot q(t)\big)\,\dot v(t)\,\mathrm dt\\
-\int_a^d \lambda_0(t)\Big[\left.\frac{\partial\theta}{\partial
q}\right\vert_{(t,q(t))}\!\!\!\!\!\big(v(t),\dot q(t)\big)+
\theta_{(t,q(t))}\,\dot v(t)\Big]\,\mathrm dt\\
-\sigma_a\Big[\left.\frac{\partial \theta}{\partial
q}\right\vert_{(a,q(a))}\!\!\!\!\!\big(v(a),\dot
q(a)\big)+\theta_{(a,q(a))}\,\dot v(a)\Big]=0.
\end{multline}
From Corollary~\ref{thm:bootstrap} and formula \eqref{eq:strong}
it follows that $\frac{\partial\widetilde L}{\partial \dot q}$ is
of class $C^1$, and we can thus use integration by parts in
\eqref{eq:standardcomp3} to obtain an equality of the form:
\begin{equation}\label{eq:trambolho}
\int_a^du(t)\,v(t)\;\mathrm dt+\sigma_a\,\theta_{(a,q(a))}\,\dot v(a)=0,
\end{equation}
for some $u\in C^0\big([a,d],{\R^n}^*\big)$, whenever $v$ is
chosen with $v(a)=0$. By considering arbitrary $v$ supported in
$\left]a,d\right[$, from \eqref{eq:trambolho} we obtain that
$u\equiv0$ in $[a,d]$, so that the integral in
\eqref{eq:trambolho} vanishes for all $v$. Now, we can choose $v$
with $v(a)=0$ and $\dot v(a)$ arbitrary, so that
\eqref{eq:trambolho} implies that $\sigma_a=0$, because
$\theta_{(a,q(a))}$ is surjective. This concludes the proof.
\end{proof}
\end{section}



\end{document}